\documentclass[12pt,a4paper]{amsart}

\usepackage{fullpage}
\usepackage{amssymb}  
\usepackage{latexsym} 
\usepackage{comment}
\usepackage{color}
\usepackage{enumerate}
\usepackage{colonequals}
\usepackage[all]{xy}

\usepackage{charter,euler} 
\DeclareSymbolFont{bchoperators}{T1}{bch}{m}{n}
\SetSymbolFont{bchoperators}{bold}{T1}{bch}{b}{n}
\makeatletter
\renewcommand{\operator@font}{\mathgroup\symbchoperators}
\makeatother

\usepackage{titlesec} 
\titleformat{\section}{\normalfont\bfseries\filcenter}{\thesection}{1em}{}
\titleformat{\subsection}{\normalfont\bfseries}{\thesubsection}{1em}{}
\titleformat{\subsubsection}{\normalfont\bfseries}{\thesubsubsection}{1em}{}

\newcommand{\Z}{{\mathbb Z}}
\newcommand{\Q}{{\mathbb Q}}
\newcommand{\F}{{\mathbb F}}
\newcommand{\BP}{{\mathbb P}}
\newcommand{\CN}{{\mathcal N}}
\newcommand{\disc}{\operatorname{disc}}
\newcommand{\Gal}{\operatorname{Gal}}
\newcommand{\eps}{\varepsilon}
\renewcommand{\sl}{\operatorname{\mathfrak{sl}}}

\newtheorem{Theorem}{Theorem}[section]
\newtheorem{Lemma}[Theorem]{Lemma}
\newtheorem{Proposition}[Theorem]{Proposition}

\newtheorem{Conjecture}[Theorem]{Conjecture}

\theoremstyle{definition}
\newtheorem{Definition}[Theorem]{Definition}

\newtheorem{Remark}[Theorem]{Remark}
\newtheorem{Question}[Theorem]{Question}

\numberwithin{equation}{section}

\definecolor{darkgreen}{rgb}{0,0.5,0}

\usepackage[
        colorlinks, citecolor=darkgreen,
        backref,
        pdfauthor={Ivan Penkov, Michael Stoll},
        pdftitle={Prime numbers and dynamics of the polynomial x²-1},
]{hyperref}
\usepackage[alphabetic,backrefs,lite]{amsrefs} 

\begin{document}

\title{Prime numbers and dynamics of the polynomial $x^2 - 1$}

\author{Ivan Penkov}
\address{(I. P.) Constructor University, 28759 Bremen, Germany}
\email{ipenkov@constructor.university}
\urladdr{\url{http://math.jacobs-university.de/penkov/}}

\author{Michael Stoll}
\address{(M. S.) Mathematisches Institut,
         Universität Bayreuth,
         95440 Bayreuth, Germany}
\email{Michael.Stoll@uni-bayreuth.de}
\urladdr{\url{http://www.mathe2.uni-bayreuth.de/stoll/}}

\date{November 11, 2025}

\keywords{Arithmetic dynamics, quadratic polynomial, prime divisors, Lie algebra}
\subjclass[2020]{37P05, 37P15; 15B30, 17B45, 11Y99}

\begin{abstract}
  Let $n \in \Z_{\geqslant 2}$. By~$P(n)$ we denote the set of all prime divisors of the integers
  in the sequence $n, n^2-1, (n^2-1)^2-1, \dots$. We ask whether the set~$P(n)$ determines~$n$
  uniquely under the assumption that $n \neq m^2-1$ for $m \in \Z_{\geqslant 2}$.
  This problem originates in the structure theory of infinite-dimensional Lie algebras. We show
  that the sets~$P(n)$ generate infinitely many equivalence classes of positive integers under
  the equivalence relation $n_1 \sim n_2 \iff P(n_1) = P(n_2)$. We also prove that the sets~$P(n)$
  separate all positive integers up to~$2^{29}$, and we provide some heuristics on why the answer
  to our question should be positive.
\end{abstract}

\maketitle


\section{Statement of the problem} 

We consider the following question.

\begin{Question} \label{Q1}
  Let $n > 1$ be an integer, not of the form $n = m^2 - 1$.
  Consider the sequence $(a_k)_{k \geqslant 0} = (a_k(n))_{k \geqslant 0}$
  defined recursively by $a_0 = n$, $a_{k+1} = a_k^2 - 1$.
  Let $P(n)$ be the set of all prime divisors of all~$a_k$.
  Is $n$ determined uniquely by $P(n)$?
\end{Question}

Note that if $p$ divides~$a_k$, then $p$ will divide infinitely many terms of the
sequence, as the sequence considered modulo~$p$ will be periodic $0 \mapsto -1 \mapsto 0$
from that point on. In particular, $P(n^2-1) = P(n)$, and so $P(a_k(n)) = P(n)$ for all~$k$.

The above question arose in the structure theory of certain infinite-dimensional Lie algebras,
see ~\cite{PH}. More precisely, for any $n \in \Z_{\geqslant 2}$ one has an interesting chain
of inclusions of Lie algebras
\[ \sl(n) \subset \sl(n^2-1) \subset \sl((n^2-1)^2-1) \subset \dots \,, \]
where the natural representation of each Lie algebra restricts to the adjoint representation
of the preceding Lie algebra. It is natural to ask when the direct limits of such chains are
isomorphic as Lie algebras. As explained in~\cite{PH}, the notion of Dynkin index allows to infer
that a positive solution to Question~\ref{Q1} implies that the direct limit Lie algebras arising
from sequences as above starting respectively with $\sl(n_1)$ and~$\sl(n_2)$ for $n_1 < n_2$ are
isomorphic if and only if $n_2$ is a member of the sequence $n_1, n_1^2-1, (n_1^2-1)^2-1, \dots$.

An equivalent formulation of Question~\ref{Q1} is as follows.
\begin{Question} \label{Q2}
  For a prime~$p$, let $\bar{S}(p)$ be the subset of the finite field $\F_p$
  consisting of all~$a$ such that iterating $x \mapsto x^2-1$ on~$a$ eventually produces~$0$.
  Denote by~$S(p)$ the preimage of~$\bar{S}(p)$ in~$\Z$.
  Is it true that for each subset~$T$ of the set of all primes the inequality
  \[ \#\bigl(\{n \in \Z_{>0} : n \neq m^2-1\} \cap \bigcap_{p \in T} S(p) \cap \bigcap_{p \notin T} (\Z \setminus S(p))\bigr) \leqslant 1 \]
  holds?
\end{Question}

On heuristic grounds (which we will detail below in Section~\ref{S:heuristics}),
this appears to be very likely, but we also expect it to be very hard to prove.

The purpose of this paper is to gather some experimental evidence and to propose some
heuristics that support a positive answer to the questions above.

We record here that some primes occur in all sets~$P(n)$ and so do not provide
useful information toward the answer of Question~\ref{Q1}.

\begin{Lemma} \label{L:no_info}
  For all~$n$, $\{2, 3, 7, 23, 19207\} \subset P(n)$, and these are the only primes
  below~$10^5$ with that property.
\end{Lemma}

\begin{proof}
  One checks by a computation that $\bar{S}(p) = \F_p$ for the five primes
    occurring in the statement and for no other primes $p < 10^5$.
\end{proof}

This begs another question:

\begin{Question} \label{Q3}
  Is $\bigcap_n P(n)$ finite or infinite?
\end{Question}

See Section~\ref{S:heuristics} for some heuristics regarding the likely answer.

\begin{Remark}
  In principle, the same questions can be asked for any polynomial with integral
  coefficients. We are not aware of any other specific polynomial for which it
  is possible to actually provide an answer to Question~\ref{Q1} with proof (other than
  trivial special cases like monomials).

  We do want to point out, however, that the polynomial $x^2 - 1$ is somewhat
  special in that its critical point~$0$ is periodic. This has the effect
  that a prime dividing one term of the iteration sequence will divide infinitely
  many of them, which gives the problem a distinctly different flavor compared
  to the generic case.
\end{Remark}

\medskip

This paper is structured as follows. In Section~\ref{S:P_infinite} we show that
all sets~$P(n)$ are infinite. Then in Section~\ref{S:finite} we provide experimental
evidence in favor of a positive answer to Questions \ref{Q1} and~\ref{Q2}.
In Section~\ref{S:infinite} we show that there are infinitely many distinct
sets~$P(n)$. In Section~\ref{S:heuristics} we present some heuristic
considerations pertaining to the questions asked above.
Finally, Section~\ref{S:speculation} has some rather speculative thoughts
on how one might try to apply a result due to Huang to reduce Question~\ref{Q1}
to Vojta's Conjecture.

\subsection*{Acknowledgments}

We used the Computer Algebra System Magma~\cite{Magma} for the computations
described in this paper. We thank the Centre for Advanced Academic Studies in Dubrovnik,
where most of this work was done. I.P.'s work was supported in part by the DFG grant PE~980/9-1.
We are grateful to Rafe Jones for information on various variants of the Basilica group
and pointers to the literature and to the anonymous referee for constructive comments
and additional references.


\section{The sets $P(n)$ are infinite (but likely sparse)} \label{S:P_infinite}

Before we state and prove the claim in the title of this section,
we need some auxiliary results.

\begin{Lemma} \label{L:prelim}
  Let $n \in \Z_{\geqslant 1}$ and let $a_k = a_k(n)$ for $k \geqslant 0$ be defined as above.
  Then for all $k \geqslant 0$ we have the divisibility $a_k \mid a_{k+2}$.
\end{Lemma}

\begin{proof}
  Note that $a_{k+2} = a_{k+1}^2 - 1 = (a_k^2 - 1)^2 - 1 = a_k^2 (a_k^2 - 2)$.
\end{proof}

We set $f \colonequals x^2 - 1 \in \Z[x]$ and write $f^m$ for its $m$-th iterate
(i.e., $f^0 = x$ and $f^{m+1} = (f^m)^2 - 1 = f^m(x^2 - 1)$)-

\begin{Lemma} \label{L:irred}
  For every $m \geqslant 0$, the polynomial $f^m - 1$ is irreducible (in~$\Z[x]$).
  Furthermore, $2$ is the only ramified prime in the splitting field of~$f^m - 1$
  (for $m \geqslant 1$).
\end{Lemma}

\begin{proof}
  The constant term of~$f^m - 1$ is~$-1$ when $m$ is even and $-2$ when $m$ is odd,
  and the constant term of $f^m(x-1) - 1$ is~$-2$ when $m$ is even and $-1$
  when $m$ is odd (we use that $f^{m+2}(0) - 1 = f^{m+1}(-1) - 1 = f^m(0) - 1$).
  The image of~$f^m$ in~$\F_2[x]$ is $x^{2^m}$ when $m$ is even and~$(x + 1)^{2^m}$
  when $m$ is odd. Both observations together imply that for each~$m$
  either $f^m - 1$ or $f^m(x-1) - 1$ is irreducible by the Eisenstein criterion at the prime~$2$.

  It is not hard to show that $\disc g(x^2 - 1) = (-4)^{\deg g} g(-1) (\disc g)^2$,
  where $g \in \Z[x]$ and $\disc g$ is the discriminant of~$g$; see, e.g.,
  \cite{Jones}*{Lemma~2.6} for a more general statement. This implies that
  the discriminant of $f^m - 1$ is a power of~$2$ up to sign, so $2$ is the only
  prime that can possibly ramify in the splitting field of $f^m - 1$. There
  are no unramified nontrivial extensions of~$\Q$, so $2$ indeed has to ramify when $m \geqslant 1$.
  (See also~\cite{BJ}*{page~222}, where the result regarding the ramification
  is stated without proof.)
\end{proof}

Now we show that $P(n)$ is infinite.

\begin{Proposition}
  Let $n \in \Z_{\geqslant 2}$. Then $P(n)$ is an infinite set of prime numbers.
\end{Proposition}

\begin{proof}
  From Lemma~\ref{L:prelim}, we can deduce that the set of primes dividing one of
  $a_0, \ldots, a_{k+1}$ is the same as the set of primes dividing~$a_k a_{k+1}$.
  Now $a_k^2 - 2$ is coprime to the odd part of $a_k a_{k+1} = a_k (a_k^2 - 1)$ and is divisible
  by~$2$ at most once. This implies that unless $a_k = 2$, $a_{k+2} = a_k^2 (a_k^2 - 2)$
  has a prime divisor not dividing $a_k a_{k+1}$. So, with one possible exception
  (which occurs only when $n = 2$), each~$a_k$ contributes at least one new prime
  to~$P(n)$. In particular, $P(n)$ must be infinite.
\end{proof}

We remark that this is a special case of the much more general Thm.~6.1
in~\cite{Jones}, which builds on results by Silverman~\cite{SilvermanBook}.
A similar result for rational functions such that $0$ is preperiodic is
given in~\cite{IngramSilverman}.

On the other hand, the sets~$P(n)$ are likely sparse in the following sense.

\begin{Conjecture} \label{Conj1}
  Let $n \in \Z_{\geqslant 2}$. Then $P(n)$ is a set of prime numbers of density zero.
\end{Conjecture}

Note that we are in the exceptional case $k = -1$ of~\cite{Jones}*{Thm.~1.2(iii)}.

The main result of~\cite{BGHKST} (in the simple form as given in the abstract there)
at least implies that $P(n)$ does not have full density (by taking $\phi(t) = t^2 - 1$
and $\alpha = 0$).

Let $G_m$ denote the Galois group of~$f^m(x) - 1$ over~$\Q$. Then Conjecture~\ref{Conj1}
would follow from the following statement.

\begin{Conjecture} \label{Conj2}
  Let $\delta_m$ be the proportion of elements $\sigma \in G_m$ such that $\sigma$
  fixes at least one root of~$f^m(x) - 1$. Then $\lim_{m \to \infty} \delta_m = 0$.
\end{Conjecture}

Note that by the Main Theorem in~\cite{ABCCF}, the corresponding group~$G_m(a)$ for
$f^m(x) - a$ is the level-$m$ quotient~$M_m$ of the `arithmetic basilica group'~$M_\infty$
for all~$m$ when $a$ is outside a `thin set'. The statement of Conjecture~\ref{Conj2}
is expected to hold in these cases (Rafe Jones, private communication). When $a = 1$, the
limit group $G_\infty = \lim\limits_{\longleftarrow} G_m$ is of infinite index in~$M_\infty$,
so this case requires additional work and is still open (with the expectation
being that Conjecture~\ref{Conj2} above should hold).


\section{Separation of numbers up to a bound} \label{S:finite}

The following definition is useful.

\begin{Definition}
  Let $X$ be a positive integer and let $P$ be a set of prime numbers.
  We say that \emph{$P$ separates the numbers up to~$X$}, if the sets $P(n) \cap P$ are pairwise
  distinct for all $n \leqslant X$ not of the form $n = m^2 - 1$.
\end{Definition}

In other words, Questions \ref{Q1} and~\ref{Q2} have a positive answer when restricted to $n \leqslant X$,
and this can be verified by only considering divisibility by primes in~$P$.

With this notion, we have the following experimental data.

\begin{Theorem} \label{T:sep}
  Write $P_{\leqslant m}$ to denote the set of prime numbers $p \leqslant m$. Then:
  \begin{enumerate}[\upshape(1)]
    \item The numbers up to~$10^1$ are separated by~$P_{\leqslant 47}$.
    \item The numbers up to~$10^2$ are separated by~$P_{\leqslant 223}$.
    \item The numbers up to~$10^3$ are separated by~$P_{\leqslant 379}$.
    \item The numbers up to~$10^4$ are separated by~$P_{\leqslant 919}$.
    \item The numbers up to~$10^5$ are separated by~$P_{\leqslant 2137}$.
    \item The numbers up to~$10^6$ are separated by~$P_{\leqslant 3001}$.
    \item The numbers up to~$10^7$ are separated by~$P_{\leqslant 4793}$.
    \item The numbers up to~$10^8$ are separated by~$P_{\leqslant 5791}$.
  \end{enumerate}
\end{Theorem}

\begin{proof}
  We run the following algorithm. The input is $X = 10^k$ with $k = 1, 2, \ldots, 8$.
  \begin{enumerate}[1.]
    \item Initialize $\CN \colonequals \bigl\{\{n \in \Z_{> 1} : n \leqslant X, \neg\exists m \colon n = m^2 - 1\}\bigr\}$,
          a set of finite sets of positive integers.
    \item Set $p \colonequals 3$.
    \item Repeat the following steps until $\CN$ is empty.
          \begin{enumerate}[a.]
            \item Replace $p$ by the next larger prime number.
            \item Compute $\bar{S}(p)$.
            \item Replace each set~$N$ in~$\CN$ by the sets
                  in the list $N \cap S(p)$, $N \setminus S(p)$ that have at least two elements.
          \end{enumerate}
    \item Return~$p$.
  \end{enumerate}
  Taking into account that the primes $2$ and~$3$ do not give information by Lemma~\ref{L:no_info},
  it is clear that this algorithm will return the minimal~$p$ such that $P_{\leqslant p}$ separates
  the numbers up to~$X$ when it terminates. The algorithm does in fact terminate
  for all $X = 10^k$ with $k \leqslant 8$ and returns the bounds given in the statement.
\end{proof}

We note that the growth of the bound on the primes that are necessary to
separate the numbers up to~$X$ is numerically consistent with a growth of
order $(\log X)^2$. The considerations in Section~\ref{S:heuristics} below would
predict $(\log X)^2 \log \log X$, which is also consistent with our numbers above
($\log \log X$ grows too slowly to allow distinguishing the two possibilities
by experimental data).

We clearly get the ``best'' effect from using the information at a prime~$p$
when the set $\bar{S}(p)$ comprises close to half the elements of~$\F_p$.
So, to get a more efficient algorithm than the one used in the proof above,
we do the following.
\begin{enumerate}[1.]
  \item Pre-compute the sets $\bar{S}(p)$ for all $p$ up to a suitable bound.
  \item Sort the list of pairs $(p, \bar{S}(p))$ by increasing value of
        $|\#\bar{S}(p)/p - 1/2|$.
  \item Use the primes in the order that is given by the sorted list of pairs.
\end{enumerate}

The effect is that we can get similar results with less computation, because
we need fewer primes to get separation.

For example, considering the primes up to~$10\,000$, the first ten primes in the
sorted list, together with the value of~$|\#\bar{S}(p)/p - 1/2|$, are
\begin{gather*}
  (2713, 0.00350), \quad (2137, 0.00726), \quad (1399, 0.0232), \quad
  (5927, 0.0534), \quad (8681, 0.0637) \\
  (4799, 0.0741), \quad (3079, 0.0746), \quad (71, 0.0775), \quad
  (919, 0.0833) \quad (7951, 0.0875) \,.
\end{gather*}

The actual splitting of the sets in~$\CN$ can be done in a breadth-first (like in
the algorithm in the proof above) or in a depth-first way. The latter is
more space-efficient, but there is no significant difference in run times
(as long as there is sufficient memory available; see below). The time
complexity should be $\asymp X \log X$ (this is corroborated by the running times),
with a memory requirement of $\asymp X$.

Using this improved algorithm, we can show:

\begin{Theorem} \label{T:2_29}
  The primes up to~$10\,000$ separate the numbers up to~$2^{29} \approx 5.37 \cdot 10^8$.
\end{Theorem}

For comparison, our
\href{https://www.mathe2.uni-bayreuth.de/stoll/magma/Penkov_question.magma}{Magma implementation}~\cite{Code}
of the algorithm described in the proof of Theorem~\ref{T:sep} takes a bit over two hours
to verify the result for $X = 10^8$ on the second author's current laptop,
while both the breadth-first and the depth-first versions of the second algorithm
take about $35$~minutes with the same bound
(but only prove the slightly weaker result that the numbers up to~$10^8$ are
separated by the primes below~$10\,000$). The computation verifying the statement
of Theorem~\ref{T:2_29} takes a bit less than four hours (using the depth-first version
and with some other tasks being executed in parallel; the breadth-first
version requires too much memory to run in reasonable time, probably caused by
the overhead incurred when working with a very large number of sets at the same time in Magma).
We expect that a low-level implementation in C could be made sufficiently (time and space)
efficient to be able to extend the bounds further, but it is perhaps not so clear
that the additional effort spent for writing, testing and debugging such an implementation
justifies the somewhat marginal improvement in experimental evidence.


\section{Infinitely many classes} \label{S:infinite}

We can define an equivalence relation on~$\Z_{>0}$ by declaring
\[ n \sim m \mathrel{\;{:}{\Longleftrightarrow}\;} P(n) = P(m) \,. \]
Then Theorem~\ref{T:2_29} shows that there are at least
$536\,847\,742 = 2^{29} - \lfloor\sqrt{2^{29} + 1}\rfloor$
distinct equivalence classes. We can in fact show more. But first we need a lemma.

\begin{Lemma} \label{L:Smin}
  Let $p \geqslant 3$ be a prime number. Then $\{-1,0,1\} \subset \bar{S}(p)$
  (as subsets of~$\F_p$) with equality if and only if $p \equiv \pm 3 \bmod 8$.
\end{Lemma}

\begin{proof}
  Since $1$ and~$-1$ both map to~$0$ under $x \mapsto x^2 - 1$, the inclusion
  $\{-1,0,1\} \subset \bar{S}(p)$ holds for all~$p$.
  The only preimage of~$-1$ is~$0$. So $\bar{S}(p)$ is larger if and only if
  there are preimages of~$1$ in~$\F_p$. But this is equivalent to $2$ being a quadratic
  residue mod~$p$, which is well known to be equivalent to $p \equiv \pm 1 \bmod 8$.
\end{proof}

\begin{Theorem}
  There are infinitely many equivalence classes under the relation ``$\sim$'' defined above.
\end{Theorem}

\begin{proof}
  We show that the sets $P(p)$, where $p \equiv \pm 3 \bmod 8$ is a prime, are pairwise distinct.
  This follows from the observation that $p \in P(p)$ and that for $q < p$ with $q \equiv \pm 3 \bmod 8$ we have
  $1 < q < p - 1$, so $q \not\equiv -1, 0, 1 \bmod p$, which means that $p \notin S(q)$
  by Lemma~\ref{L:Smin}, hence $q \notin P(p)$.
  (A similar argument works for $P(p-1)$  or $P(p+1)$ with the same primes.)
\end{proof}

One can ask whether primes $p \equiv \pm 3 \bmod 8$ are sufficient to separate
all positive integers not of the form $m^2 - 1$. However, the answer is ``no''. Indeed note that
\[ P'(n) \colonequals P(n) \cap \{p : p \equiv \pm 3 \bmod 8\} = \{p : p \equiv \pm 3 \bmod 8, p \mid (n-1) n (n+1)\} \,, \]
in particular,
\begin{align*}
   P'(2) = P'(7) = P'(17) &= \{3\} \\
   P'(4) = P'(5) = P'(6) = P'(9) = P'(16) &= \{3, 5\} \\
   P'(10) = P'(11) = P'(21) &= \{3, 5, 11\} \,.
\end{align*}
Moreover, looking at the sets~$P'(n)$ for $n$ up to~$10\,000$
seems to suggest that any given finite set of primes $p \equiv \pm 3 \bmod 8$
that contains~$3$ occurs infinitely often.


\section{Heuristics} \label{S:heuristics}

From the discussion in Section~\ref{S:finite}, it is reasonable to consider asymptotic
properties of the distribution of the relative sizes $\#\bar{S}(p)/p$ as $p$ gets large.

Here is a heuristic model for the size of~$\bar{S}(p)$:
We start with $1 \to 0 \leftrightarrow -1$ and follow~$1$ backwards.
For a given~$x \neq -1$, the chances that $x+1$ is or is not a square in~$\F_p$ are equal,
so we add two preimages with a probability of~$1/2$ and recurse.
This is related to some of the arguments in~\cite{HeathBrown},
where it is shown that the forward orbit of~$0$ under~$x^2 + c$ on a finite
field of odd size~$q$ has size $\ll q/\log\log q$.

We write
\[ F(z) = \sum_{n \ge 0} \BP(\#\bar{S}(p) = n) z^n = z^2 G(z) \in \Q[\![z]\!] \,, \]
where
\[ G(z) = \sum_{n \ge 0} g_n z^n \]
with $g_n$ the probability that a binary tree that is randomly generated by the
following procedure has $n$ nodes in total: start with a root node and for each
node in the tree, add two children with probability~$1/2$. This tree is obtained
by removing $-1$ and~$0$ from~$\bar{S}(p)$; $x$ is a child of $y$ when $y = x^2 - 1$.

Alternatively, we obtain the same distribution starting from the root node
by leaving it alone with probability~$1/2$ and otherwise adding two random trees
of the same kind to it as the left and right sub-trees. This leads to the equation
\[ G(z) = \frac{z}{2} \bigl(1 + G(z)^2\bigr)  \quad\Longrightarrow\quad  G(z) = \frac{1 - \sqrt{1-z^2}}{z} \,, \]
so
\[ F(z) = z \bigl(1 - \sqrt{1 - z^2}\bigr) = \sum_{n \geqslant 1} (-1)^{n+1} \binom{1/2}{n} z^{2n+1} \,. \]
The coefficient of~$z^{2n+1}$ is
\[ \left|\binom{1/2}{n}\right| = \frac{1}{2n-1} \left|\binom{-1/2}{n}\right|
     = \frac{4^{-n}}{2n-1} \binom{2n}{n}
     \sim \frac{1}{2n \sqrt{\pi n}} \,.
\]

For the regime of relative sizes $< 1 - \varepsilon$, the heuristic model above
should be fairly accurate. Indeed, counting the primes $p < 10^5$ such that $\#\bar{S}(p) = 2n+1$
shows a reasonably good agreement with the prediction from the model for, say, $n \leqslant 20$
(for larger~$n$ the numbers of primes are too small for a meaningful comparison).

The model would predict that for $0 < a < b < 1$, we should expect
there to be about
\[ \operatorname{const} \int_3^X \left(\int_{ax}^{bx} \frac{dt}{t^{3/2}}\right)\, d\pi(x)
     \approx \operatorname{const} \frac{\sqrt{X}}{\log X} \left(\frac{1}{\sqrt{a}} - \frac{1}{\sqrt{b}}\right)
\]
primes~$p \leqslant X$ such that $a p \leqslant \#\bar{S}(p) \leqslant b p$.

For instance, we expect there to be infinitely many primes~$p$ such that
\[ \left|\frac{\#\bar{S}(p)}{p} - \frac{1}{2}\right| \leqslant \frac{1}{10} \,. \]
The first few primes satisfying this inequality are
\begin{gather*}
  5, \; 71, \; 919, \; 1399, \; 2137, \; 2713, \; 3079, \; 4799, \; 5927, \; 7951, \; 8681, \; 10271, \; 10711, \; 11369, \\
  12487, \; 12577, \; 22409, \; 22871, \; 24623, \; 24631, \; 27647, \; 29641, \; 46457, \; 54751, \; 84559, \\
  87583, \; 99929, \; 103703, \; 105449, \; 106753, \; 120199, \; 120607, \; 123289, \; 131111, \; 147703\,.
\end{gather*}
This fits reasonably well with the heuristic growth.

Each such prime will lead to a nearly $50$-$50$ split of the numbers~$n$, and so
the very likely fact that there are infinitely many of them gives another strong indication
that the sets~$P(n)$ do actually separate all positive integers~$n$ not of the form~$m^2 - 1$.


\medskip

For the ``probability'' that $\#\bar{S}(p) = p$ (i.e., $\bar{S}(p) = \F_p$), the model
predicts something like~$O(p^{-3/2})$, which would suggest that the answer to Question~\ref{Q3}
is ``finite''.

However, the model does not take into account that $\#\bar{S}(p) \leqslant p$. If we
include the cases where the size would be larger than~$p$ according to the model,
the ``probability'' grows to~$p^{-1/2}$, which would indicate that the answer is ``infinite''!
A refined model will be necessary to obtain reasonable predictions in the regime
of $\bar{S}(p)$ close to maximal. So we now develop
another heuristic regarding the existence of infinitely many primes with $\bar{S}(p) = \F_p$.

The directed graph on the vertex set~$\F_p$ with edges $x \to y$ when $y = x^2 - 1$
consists of a number~$N(p)$ of connected components, each of which contains
precisely one cycle. Define the polynomials~$\Phi_n$ so that
\[ \prod_{d|n} \Phi_d = f^n - x \,, \]
where $f = x^2 - 1$ and $f^n$ denotes the $n$-th iterate of~$f$. So, e.g.,
\begin{align*}
  \Phi_1 &= x^2 - x - 1 \\
  \Phi_2 &= x^2 + x = x   (x + 1) \\
  \Phi_3 &= x^6 + x^5 - 2 x^4 - x^3 + x^2 + 1 \\
  \Phi_4 &= x^{12} - 6 x^{10} + x^9 + 12 x^8 - 4 x^7 - 7 x^6 + 4 x^5 - 4 x^4 + x^3 + 4 x^2 - 2 x + 1 \\
  \text{etc.}
\end{align*}
At least up to $m = 6$, we have
\[ \Gal\bigl(\prod_{n \leqslant m} \Phi_n\bigr) = \prod_{n \leqslant m} \Gal(\Phi_n) \]
(where $\Gal(g)$ denotes the Galois group of a polynomial~$g$ over~$\Q$)
and (for $n \neq 2$)
\[ \Gal(\Phi_n) = G_n \colonequals C_n \wr S_{k_n} \quad \text{with} \quad
      k_n = \frac{\deg \Phi_n}{n} = \frac{1}{n} \sum_{d \mid n} \mu\Bigl(\frac{n}{d}\Bigr) 2^d \,.
\]
Here, $C_n \wr S_{k_n}$ denotes the wreath product, i.e., the semidirect product of $C_n^{\#S_{k_n}}$
with $S_{k_n}$ acting via permutation of the factors, and $C_n$ is the cyclic group of order~$n$.

If we assume that this remains true for larger~$m$ (this is the case for generic polynomials~$f$),
then the density of primes~$p$
such that $f$ has no cycle of length $n$ ($\neq 2$) in~$\F_p$ is (by the Chebotarev Density Theorem)
the fraction of elements in~$G_n$ without a fixed point, which is
\[ \sum_{k=0}^{k_n} \left(\sum_{j=0}^{k_n-k} \frac{(-1)^j}{j!}\right) \frac{(1 - 1/n)^k}{k!} \,. \]
This is quite close to
\[ e^{-1} e^{1 - 1/n} = e^{-1/n} \,. \]
So under this model and assuming that the events ``$f$ has a cycle of length~$m$ on~$\F_p$''
are independent for $3 \leqslant m \leqslant p$,
the expected ``probability'' for a prime~$p$ to have $N(p) = 1$ is close to
\[ \frac{1}{2} \exp\left(-\sum_{3 \leqslant n \leqslant p} \frac{1}{n}\right)
    \sim \frac{1}{2} \exp\bigl(-(\log p + \gamma - 3/2)\bigr)
    = \frac{c}{p}
\]
for a constant~$c$. This would lead us to expect roughly ~$c \log \log X$  such primes
up to~$X$. In particular, this indicates that the answer to Question~\ref{Q3} is
``infinite''.

A bit more realistically, since (except for~$-1$) every element in a cycle
has \emph{two} preimages, we can stop at $m = (p-3)/2$  (subtract $3$ for $1, 0, -1$),
which basically doubles~$c$.
We'd then actually expect a prime $\approx 163$ and one $\approx 2130$ and then one
around~$100\,000$ in the intersection of all sets~$P(n)$. Of course, these are not
precise predictions, and they serve only as an indication of the expected growth
of the numbers in~$\bigcap_n P(n)$.

See also~\cite{JKMT} for some general results on the proportion of preperiodic points
mod~$p$ under iteration of polynomials or rational functions.

\begin{Remark}
  The use of heuristics based on the assumption that polynomial maps
  on finite fields behave randomly is fairly standard in the study of
  arithmetic dynamical systems; see, e.g., \cite{BGHKST}*{Section~5}.
\end{Remark}


\section{A speculation} \label{S:speculation}

Assuming Vojta's Conjecture, Huang~\cite{Huang}, inspired by Silverman's
re-interpretation in~\cite{Silverman} of a result by Bugeaud, Corvaja
and Zannier~\cite{BCZ}, shows a general result that says the following.
Let $f_1, f_2 \in \Z[x]$ be polynomials of degree~$d$ satisfying some mild condition,
let $\eps$ be a fixed positive real number, and let $a$ and~$b$ be integers.
If the sequence $\bigl((f_1^n(a), f_2^n(b)\bigr)_{n \ge 0}$ is Zariski-dense
in the affine plane, then
\[ \gcd\bigl(f_1^n(a), f_2^n(b)\bigr) \ll e^{\eps d^n} \,. \]
(We again write $f_j^n$ for the $n$th iterate of~$f_j$.)

We can try to apply this with $f_1 = f_2 = f^2$ (where $f^2(x) = x^2 (x^2 - 2)$
is the second iterate of $f = x^2 - 1$); it is not hard to see that a positive
answer to the analogue of Question~\ref{Q1} for~$f^2$ implies a positive answer
to the original question.

Write $P^2(n)$ for the set of primes dividing $a_{2k}(n)$ for some~$k \ge 0$.
Assume that $P^2(n) = P^2(n')$ for $n, n' \ge 2$. \emph{If} we could show that this
implies that there is some $\delta > 0$ such that
\[ \gcd\bigl(f^{2k}(n), f^{2k}(n')\bigr) > e^{\delta 4^k} \]
for infinitely many~$k$, then Huang's result would show (if we also assume
Vojta's Conjecture) that the sequence $\bigl((f^{2k}(n), f^{2k}(n'))\bigr)_{k \ge 0}$
is not Zariski dense, so there is some polynomial~$F \in \Z[u, v]$ such that
\[ F\bigl(f^{2k}(n), f^{2k}(n')\bigr) = 0 \]
for all $k \ge 0$. It seems rather plausible that this would imply that $n'$
is in the $f^2$-orbit of~$n$ or vice versa.

Note that when $p \in P^2(n) \cap P^2(n')$, then there is some $k_0(p)$ such that
\[ v_p\bigl(\gcd\bigl(f^{2k}(n), f^{2k}(n')\bigr)\bigr) \ge 2^{k-k_0(p)}
    \qquad \text{for all $k \ge k_0(p)$.}
\]
Here $v_p(m)$ denotes the exponent of~$p$ in the prime factorization of~$m$.
The desired lower bound would follow if we could prove something like
\[ \sum_{p \colon k_0(p) \le k} 2^{-k_0(p)} \log p \gg 2^k \qquad \text{as $k \to \infty$} \,, \]
where the sum is over all $p \in P^2(n) = P^2(n')$ such that $k_0(p) \le k$.
This probably comes down to showing that the first index~$k$ such that
a prime~$p$ divides $a_{2k}(n)$ is not too far away (in most cases) from
the corresponding index for~$n'$. Unfortunately, this still seems to be
rather difficult.


\end{document}